\documentclass[12pt]{article}
\usepackage{amsfonts}

\newtheorem{lemma}{Lemma}
\newtheorem{theorem}{Theorem}

\begin{document}

\title{The Novikov--Veselov hierarchy of equations and integrable deformations of minimal Lagrangian tori in
${\mathbb C}P^2$}
\author{A. E. Mironov}
\date{}

\maketitle
\begin{abstract}
We associate a periodic two-dimensional Schr\"odinger operator to
every Lagrangian torus in ${\mathbb C}P^2$ and define the spectral
curve of a torus as the Floquet spectrum of this operator on the
zero energy level. In this event minimal Lagrangian tori
correspond to potential operators. We show that Novikov--Veselov
hierarchy of equations induces integrable deformations of minimal
Lagrangian torus in ${\mathbb C}P^2$  preserving the spectral
curve. We also show that the highest flows on the space of smooth
periodic solutions of the Tziz\'eica equation are given by the
Novikov--Veselov hierarchy.

\end{abstract}

\section{Introduction}

The surface $\Sigma$ in ${\mathbb C}P^2$ is called Lagrangian if
the restriction of the Fubini-Studi form on $\Sigma$ is equal to
zero. Let $S^5$ be a unit sphere in ${\mathbb C}^3$, and
 ${\cal H}:S^5\rightarrow{\mathbb C}P^2$
--- a Hopf bundle. Define a conformal Lagrangian immersion
$\varphi:\Omega\rightarrow{\mathbb C}P^2$ of the domain $\Omega\subset{\mathbb R}^2$
as a composition $r:\Omega\rightarrow S^5$ and ${\cal H}$.

The following lemma holds:

\begin{lemma}
The components $r_j$ of the vector function $r$ satisfy the Schr\"odinger equation
$$
 Lr_j=\partial_x^2r_j+\partial_y^2r_j+i(\beta_x\partial_xr_j+\beta_y\partial_yr_j)+4e^vr_j=0
$$
where $2e^v(dx^2+dy^2)$ is an induced metric on the surface $\varphi(\Omega)$,
and $\beta(x,y)$ --- a Lagrangian angle defined by the equality
$$
e^{i\beta}=dz_1\wedge dz_2\wedge dz_3(\sigma),
$$
$z_1,z_2,z_3$ --- coordinates in ${\mathbb C}^3$, $x,y$ --- coordinates in $\Omega$, $\sigma$ ---
a frame formed by the vectors $r,\frac{r_x}{|r_x|},\frac{r_y}{|r_y|}$.
\end{lemma}

Lemma 1 allows the following definition.

A Lagrangian torus defined by a doubly periodic conformal mapping
$$
 \varphi={\cal H}\circ r:{\mathbb R}^2\rightarrow{\mathbb C}P^2
$$
is called {\sl finite-gap} if the corresponding Schr\"odinger operator $L$ with periodic coefficients
is finite-gap on a zero energy level, i.e. if the Bloch functions (common eigenfunctions of $L$
and of the translation operators) of the operator $L$ on a zero energy level are parametrized
by a Riemannian surface $\Gamma$ of finite genus. The Riemannian surface $\Gamma$ is called
the {\sl spectrum} of a Lagrangian torus, and its genus --- the torus' {\sl spectral genus}.

The concept of a spectrum and of the torus' spectral genus was first introduced by Taimanov [2]
for an arbitrary smooth torus in ${\mathbb R}^3$. Here, a Dirac operator serves as an analogue to the Schr\"odinger operator.

Since the mapping $\varphi$ is doubly periodic, the components of the vector function
$r$ are Bloch functions of the operator $L$.

Finite-gap Schr\"odinger operators with respect to one energy level were first introduced by
Dubrovin, Krichever and Novikov [3]. The authors of [3] also indicate data of an inverse problem.
These are used to recover the potential and the magnetic field, and to deduce explicit formulae.

The Lagrangian angle of minimum Lagrangian surfaces in ${\mathbb C}P^2$
is constant (see, for example, [4]);
thus the potential Schr\"odinger operators
$$
 L=\partial_x^2+\partial_y^2+4e^v
$$
correspond to them.
In this case Lemma 1 is analogous to the assertion that the components of the vector functions that define a
conformal minimal immersion of a plane
domain into the three--dimensional Euclidean space are harmonic functions.

Clifford's torus in ${\mathbb C}P^2$ has the spectral genus 0.
Castro and Urbano [5] and Joyce [6] have constructed examples of minimal Lagrangian tori of spectral genera 2 and 4.
Note that the spectral genus of a finite-gap minimal Lagrangian torus is even because the spectral curve
of a finite-gap potential Schr\"odinger operator
has a holomorphic involution [1].

Sharipov [7] proved that the metric of a minimal torus in $S^5$ satisfies the Tziz\'eica equation.
He also constructed finite-gap solutions of this equation. Actually, as already noted in [8],
the Sharipov construction is suitable for the construction of all minimal Lagrangian tori in
${\mathbb C}P^2$. For this it is necessary to apply the Hopf mapping to the
mapping from ${\mathbb R}^2$ into $S^5$ that was constructed by Sharipov.

Since work [7] does not discuss the problem of periodicity of the
constructed mappings, it leaves unclear the question whether there
exist minimal Lagran\-gian tori of arbitrary spectral genus. This
gap was filled by a work of Carberry and McIntosh [9] which proved
that for any spectral genus $g$ there exists a
$(\frac{g}{2}-2)$-dimensional set of minimal Lagrangian tori in
${\mathbb C}P^2$.

By the methods used in [10] it can be shown that all minimal smoothly immersed
Lagrangian tori are finite-gap. In fact, the metric of a minimal Lagrangian torus fulfills the Tziz\'eica equation
$$
\partial_x^2v+\partial_y^2v=4e^{-2v}-4e^v.
$$
Mikhailov [11] proved that the Tziz\'eica equation is integrable
(as an infinite-dimensional Hamiltonian system). Therefore, all
its smooth real doubly perio\-dic solutions are finite-gap.
Finite-gap solutions are characterized by stationa\-rity with
respect to some higher flow [12]. In our case this follows from
the fact that the function $\partial_{t_i}v$, where $t_i$ is the
highest time, satisfies the elliptic equation
$$
 (\partial_x^2+\partial_y^2+8e^{-2v}+4e^v)\partial_{t_i}v=0
$$
on the torus ${\mathbb R}^2/\Lambda$, where $\Lambda$ is a lattice of periods.
Since the spectrum of the elliptic operator on a torus is discrete, the functions $\partial_{t_i}v$
are linearly dependent and there exists a peak time with respect to which $v$ is stationary.

The work's basic result consists of the following.
Let the map
$$
 r:{\mathbb R}^2\rightarrow S^5
$$
define a finite-gap minimal Lagrangian torus $T\subset {\mathbb C}P^2$ of spectral genus $g>4$.
Then the following is true:

\begin{theorem}
There is a mapping $\widetilde{r}(t),t=(t_1,t_2,\dots), \widetilde{r}(0)=r$,
defining a deformation of torus $T$ in the class of minimal Lagrangian tori in ${\mathbb C}P^2$.
The map $\widetilde{r}$ satisfies the equations
$$
 L\widetilde{r}=
 \partial_x^2\widetilde{r}+\partial_y^2\widetilde{r}+4e^{\widetilde v}\widetilde{r}=0,
$$
$$
 \partial_{t_n}\widetilde{r}=A_n\widetilde{r},\
$$
where $A_n$ are operators of order $(2n+1)$ on the variables $(x,y)$.
Deform the potential
$\widetilde{V}=4e^{\widetilde{v}},\widetilde{v}(0)=v,$
according to the  Novikov--Veselov hierarchy
$$
 \frac{\partial L}{\partial t_n}=[L,A_n]+B_nL,\
$$
where $B_n$ are operators of order $(2n-1)$ on the variables
$(x,y)$. The deforma\-tions $\widetilde{r}(t)$ preserve the
spectrum of torus $T$ and its conformal type.
\end{theorem}

Thus, the highest flows on the space of smooth periodic solutions of the Tziz\'eica equation are given
by the Novikov--Veselov hierarchy.

It can be shown that the deformations corresponding to the
Novikov--Veselov hierarchy's first equation leave the torus
geometrically unchanged on it's spot, while from the second
equation the deformations nontrivial.

We suggest, using the first eguation of the  Novikov--Veselov hierarchy one can construct deformations of any
arbitrary Lagrangian torus, with the minimal tori being immovable.

As was shown by Taimanov [14], the local deformations of surfaces
in ${\mathbb R}^3$ introduced in [13], under the action of the
modified Novikov--Veselov equation transform tori into tori
preserving the Wilmore functional. Distinct from our construction,
[13] defines the deformation of the tori not by the deformation of
a radius--vector, but of a Gauss mapping. The proof that the
surface remains closed under the action of such deformations
substantially uses the characteristics of the modified
Veselov--Novikov equation. However, in our case the closure of the
surfaces follows from the explicit form $\widetilde{r}(t)$ (see
below).

The proof of Theorem 1 is based on Lemma 1 and the Sharipov construction [6].

The author thanks I.A. Taimanov for valuable discussions.

\section{Proof of Theorem 1}
Since the map $\varphi$ is Lagrangian and conformal, it is easily verified [4] that
$$
 <r,r_x>=<r,r_y>=<r_x,r_y>=0,\ |r_x|^2=|r_y|^2=2e^{v},
$$
where $<.,.>$ is the Hermitian product in ${\mathbb C}^3$. Thus, from the definition of the Lagrangian angle  $\beta$
we obtain
$$
 R=
 \left(
 \begin{array}{c}
   r\\
  e^{-i\frac{\beta}{2}} \frac{r_x}{|r_x|} \\
   e^{-i\frac{\beta}{2}} \frac{r_y}{|r_y|} \\
 \end{array}\right)=
  \left(
  \begin{array}{ccc}
    r^1 & r^2 &  r^3\\
   \frac{1}{\sqrt{2}}e^{-\frac{v}{2}-i\frac{\beta}{2}}r^1_x
   &  \frac{1}{\sqrt{2}}e^{-\frac{v}{2}-i\frac{\beta}{2}}r^2_x &
    \frac{1}{\sqrt{2}}e^{-\frac{v}{2}-i\frac{\beta}{2}}r^3_x \\
    \frac{1}{\sqrt{2}}e^{-\frac{v}{2}-i\frac{\beta}{2}}r^1_y
    &  \frac{1}{\sqrt{2}}e^{-\frac{v}{2}-i\frac{\beta}{2}}r^2_y
    &  \frac{1}{\sqrt{2}}e^{-\frac{v}{2}-i\frac{\beta}{2}}r^3_y \\
  \end{array}\right)\in{\rm SU(3)},
$$
where $r^1,r^2$ and $r^3$ are components of the vector $r$.
The matrix $R$ satisfies the equations
$$
 R_x=AR,\ R_y=BR,\eqno{(1)}
$$
where matrices $A$ and $B$ have the form
$$
 A=
  \left(
  \begin{array}{ccc}
   0 & \sqrt{2}e^{\frac{v}{2}+i\frac{\beta}{2}} & 0\\
  -\sqrt{2}e^{\frac{v}{2}-i\frac{\beta}{2}} & if &
  -\frac{v_y}{2}+i(h+\frac{\beta_y}{2}) \\
  0 & \frac{v_y}{2}+i(h+\frac{\beta_y}{2}) & -if\\
  \end{array}\right)\in{\rm su(3)},
$$
$$
 B=
  \left(
  \begin{array}{ccc}
   0 & 0 & \sqrt{2}e^{\frac{v}{2}+i\frac{\beta}{2}} \\
  0 & ih & \frac{v_x}{2}+i(-f+\frac{\beta_x}{2})  \\
  -\sqrt{2}e^{\frac{v}{2}-i\frac{\beta}{2}} &
  -\frac{v_x}{2}+i(-f+\frac{\beta_x}{2}) & -ih\\
  \end{array}\right)\in{\rm su(3)},
$$
$f(x,y)$ and $h(x,y)$ are some functions.
From the zero curvature equation
$$
 A_y-B_x+[A,B]=0
$$
 follows the next lemma (see [15])
\begin{lemma}
The following equations hold:
$$
2{\cal G}_y+2{\cal F}_x=(\beta_{xx}-\beta_{yy})e^v,
$$
$$
2{\cal F}_y-2{\cal G}_x=(\beta_yv_x+\beta_xv_y)e^v,
$$
$$
\Delta v=4({\cal F}^2+{\cal G}^2)e^{-2v}-4e^{v}-2({\cal F}\beta_x+{\cal G}\beta_y)e^{-v},
$$
where ${\cal F}=fe^{v},{\cal G}=he^{v}.$
\end{lemma}

From (1) we obtain the equalities
$$
 r_{xx}=\frac{1}{2}(-4e^{v}r+r_x(2if+v_x+i\beta_x)+r_y(2ih-v_y+i\beta_y)),
$$
$$
 r_{yy}=\frac{1}{2}(-4e^{v}r+r_x(-2if-v_x+i\beta_x)+r_y(-2ih+v_y+i\beta_y)).
$$
From these equalities follows Lemma 1.

Below we consider minimal Lagrangian tori.
From Lemma 2 we obtain $\Delta {\cal F}=\Delta {\cal G}=0$; consequently,
since functions ${\cal F}$ and ${\cal G}$ are doubly periodic, ${\cal F}$ and ${\cal G}$ are constants
and from Lemma 2 follows the Tziz\'eica equation.

Consider the following equations with the spectral parameter $\lambda$
$$
 \partial_zR(\lambda)=A(\lambda)R(\lambda),\
 \partial_{\bar{z}}R(\lambda)=B(\lambda)R(\lambda),\eqno{(2)}
$$
where $z=x+iy$,
$$
 A(\lambda)=
  \left(
  \begin{array}{ccc}
   0 & 1 & 0\\
  0 & v_z & -\frac{i}{\lambda}{e^{-v}} \\
  -e^{v} & 0 & 0\\
  \end{array}\right),\
 B(\lambda)=
  \left(
  \begin{array}{ccc}
   0 & 0 & 1 \\
  -e^{v} & 0 & 0 \\
  0 &
  -i\lambda e^{-v} & v_{\bar{z}}\\
  \end{array}\right),\
$$
$$
 R(\lambda)=\left(
  \begin{array}{c}
   r(\lambda) \\
   r_z(\lambda) \\
   r_{\bar{z}}(\lambda)\\
  \end{array}\right).
$$
For $\lambda=1$  equations (2) are equivalent to the equations (1) ($\beta=0$),
and the zero curvature equation for matrices $A(\lambda)$ and $B(\lambda)$ are equivalent
to the Tziz\'eika equation for any $\lambda$.
In the case of the finite gap solutions of the equations (2) there is a matrix $W(x,y,\lambda)$
rationally depending on $\lambda$ [16] such that
$$
 W_z=[A(\lambda),W],\ W_{\bar{z}}=[B(\lambda),W].
$$
The coefficients of the rational function on $\lambda$ and $\mu$
$$
 Q(\lambda,\mu)={\rm det}(W-\mu E),
$$
where $E$ is a unit matrice, do not depend on $x$ and $y$.
The spectrum of a minimal Lagrangian torus is given in the $(\lambda,\mu)$--plane
by the equation $Q(\lambda,\mu)=0$. Hence, the spectrum is the three--sheeted cover of the
$\lambda$-plane, i.e. the spectrum is a trigonal curve.

For the construction of the minimal finite-gap Lagrangian tori, recall the following construction.
Finite-gap real potential Schr\"odinger operators are built on the follwing spectral data [1]:
$\Gamma$ is a nonsingular Riemannian surface of even genus $g=2g_0$, two marked points
$\infty_1, \infty_2\in\Gamma$, a nonspecial divisor $D=P_1+\dots+P_g,$
local parameters $k_1^{-1}$ and $k_2^{-1}$ near points $\infty_1$ and $\infty_2.$
Surface $\Gamma$ should have a holomorphic involution
$$
 \sigma:\Gamma\rightarrow\Gamma,\ \sigma^2=1,
$$
with two fixed points
$\infty_1$ and $\infty_2$ such that
$$
 \sigma(k_s^{-1})=-k_s^{-1}, \
 D+\sigma D=\infty_1+\infty_2+K,
$$
where $s=1,2, \ K$ is a canonical class on $\Gamma$.
In order for $L$ to be real, the surface $\Gamma$ must have an antiholomorphic involution commutative to $\sigma$
$$
 \tau:\Gamma\rightarrow\Gamma,\ \tau^2=1,
$$
such that
$$
 D=\tau(D),\ \tau(\infty_1)=\infty_2,\ k_1(\tau(P))=\overline{k_2(P)}.
$$
There is a unique function $\psi(P,x,y)$ called Baker--Akhiezer function which is
meromorphic on
$\Gamma\backslash\{\infty_1,\infty_2\}$ and has simple poles on divisor $D$
and the following asymptotics
$$
 \psi(P,x,y)=\exp(k_1z)\left(1+\frac{\xi(x,y)}{k_1}+\dots\right),\
 P\rightarrow\infty_1,
$$
$$
 \psi(P,x,y)=\exp(k_2\bar{z})\left(1+\frac{\eta(x,y)}{k_2}+\dots\right),\
 P\rightarrow\infty_2.
$$
Function $\psi$ satisfies the Schr\"odinger equation
$$
 \partial_x^2\psi+\partial_y^2\psi+4e^v\psi=0
$$
where  $e^v=-\xi_{\bar{z}}=-\eta_z.$

Below we explain the Sharipov construction to build finite-gap
solutions of the Tziz\'eica equation [7] (in [7] instead of the
anti-holomorphic involution $\tau$ we consider, in our
terminology, an anti-holomorphic involution $\sigma\tau$). Let
curve $\Gamma$ have a meromorphic function $\lambda$ with the
divisor of zeros and poles $3\infty_1-3\infty_2$ such that
$$
 \lambda(\sigma(P))=-\lambda(P),\ \lambda(\tau(P))\overline{\lambda(\sigma(P))}=1.\eqno{(3)}
$$
Choose $k_1$ and $k_2$ such that in the vicinity of $\infty_1$ and $\infty_2$
function $\lambda$ has the form
$$
 \lambda=ik_1^{-3}, \ P\rightarrow\infty_1,\
$$
$$
  \lambda=\frac{k_2^3}{i},\
 P\rightarrow\infty_2.
$$
The choice of such spectral data provides for the smoothness and
realness of the potential of the Schr\"odinger operator. From the
uniqueness of the Baker--Akhiezer function follow the equalities
$$
 \psi_{zz}=\frac{\xi_{z\bar{z}}}{\xi_{\bar{z}}}\psi_z+\frac{k_1^3}{\xi_{\bar{z}}}
 \psi_{\bar{z}}=v_z\psi_z-\frac{i}{\lambda}e^{-v}\psi_{\bar{z}},\eqno{(4)}
$$
$$
 \psi_{\bar{z}\bar{z}}=\frac{k_2^3}{\eta_z}\psi_z
 +\frac{\eta_{z\bar{z}}}{\eta_z}\psi_{\bar{z}}=
 -e^{-v}i\lambda\psi_z+v_{\bar{z}}\psi_{\bar{z}},\eqno{(5)}
$$
$$
\psi_{z\bar{z}}=\xi_{\bar{z}}\psi=\eta_z\psi=-e^v\psi,\eqno{(6)}
$$
$$
 \psi(P)=\overline{\psi(\tau(P))},\ \psi_{\bar{z}}(P)=\overline{\psi_z(\tau(P))},\
 \psi_z(P)=\overline{\psi_{\bar{z}}(\tau(P))}.\eqno{(7)}
$$
Consider the function
$$
 F(P,Q)=<e(P),e(Q)>,
$$
where
$e(P)=(\psi(P),\psi_z(P)e^{-\frac{v}{2}},\psi_{\bar{z}}(P)e^{-\frac{v}{2}})$.
From (7) obtain
$$
 F(P,Q)=\psi(P)\psi(\tau(Q))+\psi_z(P)\psi_{\bar{z}}(\tau(Q))e^{-v}
 +\psi_{\bar{z}}(P)\psi_z(\tau(Q))e^{-v}.
$$
From (3)--(6) obtain
$$
 F_z(P,Q)=-ie^{-2v}\left(\frac{1-\lambda(P)\overline{\lambda(Q)}}{\lambda(P)}\right)
 \psi_{\bar{z}}(P)\psi_{\bar{z}}(\tau(Q)),
$$
$$
 F_{\bar{z}}(P,Q)=-ie^{-2v}
 \left(\frac{\lambda(P)\overline{\lambda(Q)}-1}{\overline{\lambda(Q)}}\right)
 \psi_z(P)\psi_z(\tau(Q)).
$$
Function $\lambda$ gives a three--sheeted cover of ${\mathbb
C}P^1$ by the curve $\Gamma$. Let
$\lambda(P_1)=\lambda(P_2)=\lambda(P_3)=1$. Then the function
$F(P_i,P_j)$ does not depend on $x$ and $y$.

Put
$$
 r_j=C_j\psi(P_j),
$$
where $C_j=\frac{1}{|\psi(P_j)|}$. In this case the equations (1)
are fulfilled, where matrices $A$ and $B$ belong to the Lie
algebra ${\rm su(3)}$ (see [7]).

By the following asymptotics, the Baker--Akhiezer function defines the integrable deformations of torus $T$:
$$
 \psi(P,x,y,t)=\exp(k_1z+\sum_{n=1}^{\infty}k_1^{2n+1}t_n')
 \left(1+\frac{\xi(x,y,t)}{k_1}+\dots\right),\
 P\rightarrow\infty_1,
$$
$$
 \psi(P,x,y,t)=\exp(k_2\bar{z}+\sum_{n=1}^{\infty}k_2^{2n+1}\bar{t}_n')
 \left(1+\frac{\eta(x,y,t)}{k_2}+\dots\right),\
 P\rightarrow\infty_2,
$$
where $t'_j=t_j+it_j.$ The function $\psi(P,x,y,t)$ has the
properties (4)--(7).

The function $\psi$ can be extracted in terms of a Prym theta-function
of the involution $\sigma$ (see [1]).
There is a basis of cycles $a_1,\dots,a_g,b_1,\dots,b_g$ on $\Gamma$ such that
$$
 \sigma(a_i)=a_{i+g_0},\ \sigma(b_i)=-b_{i+g_0},i=1,\dots,g_0
$$
and a corresponding basis of Abelian differentials
$\omega_1,\dots,\omega_g$, with the properties
$$
\int_{a_j}\omega_k=2\pi i\delta_{jk}.
$$
The Prym variety of $(\Gamma,\sigma)$ is
$$
 P={\mathbb C}^{g_0}/\{2\pi i{\mathbb Z}^{g_0}+\Omega{\mathbb Z}^{g_0}\},
$$
where the components of the symmetric matrix $\Omega$ are periods
of the following differentials: $\Omega_{ij}=\int_{b_j}\eta_i, \
\eta_i=\omega_i+\omega_{i+g_0}$. Let $\eta(P)$ denote the map
$$
 \eta:\Gamma\rightarrow P,\ \eta(P)=\left(\int_{P_0}^P\eta_1,\dots,\int_{P_0}^P\eta_{g_0}\right),
$$
where $P_0\in\Gamma$ is some fixed point. By $\Omega_k$ and
$\widetilde{\Omega}_k,\ k=0,1,\dots$
we denote meromorphic differentials on $\Gamma$ with unique poles in
$\infty_1$ and $\infty_2$ of the form $d(k_s^{-(2k+1)}), s=1,2$ and normalized by
$\int_{a_j}\Omega_k=\int_{a_j}\widetilde{\Omega}_k=0.$ Let
$$
 V_k=\left(\int_{b_1}\Omega_k,\dots,\int_{b_{g_0}}\Omega_k\right),\
 \widetilde{V}_k=\left(\int_{b_1}\widetilde{\Omega}_k,\dots,\int_{b_{g_0}}
 \widetilde{\Omega}_k\right).
$$
The theta-function of the Prym variety is defined by the convergent series
$$
 \theta(z)=\sum_{n\in{\mathbb Z}^{g_0}}\exp\left(\frac{1}{2}<\Omega n,n>+<z,n>\right),
$$
$z=(z_1,\dots,z_{g_0})\in{\mathbb C}^{g_0}.$
The theta-function has the properties of periodicity
$$
 \theta(z+2\pi in+\Omega m)=\exp\left(-\frac{1}{2}<\Omega m,m>+<z,m>\right)\theta(z),
$$
$m,n\in{\mathbb Z}^{g_0}$.
The function $\psi$ has the following form (see [1]):
$$
 \psi=\frac{\theta(\eta(P)+zV_0+\bar{z}\widetilde{V}_0+t'_1V_1+\bar{t'_1}\widetilde{V}_1+\dots-e)}
{\theta(\eta(P)-e)\theta(zV_0+\bar{z}\widetilde{V}_0+t'_1V_1+\bar{t'_1}\widetilde{V}_1+\dots-e)}
$$
$$
 \times
 \exp\left(z\left(\int_{P_0}^P\Omega_0-\alpha_0\right)+
 \bar{z}\int_{\infty_1}^{P}\widetilde{\Omega}_0+
 t'_1\left(\int_{P_0}^P\Omega_1-\alpha_1\right)
 +\bar{t'}_1\int_{\infty_1}^{P}\widetilde{\Omega}_1+\dots\right),
$$
$\alpha_j$ are some constants, $e\in{\mathbb C}^{g_0}$ is some vector.
Put
$$
 \widetilde{r}_j(t)=C_j(t)\psi(P_j,x,y,t),
$$
where $C_j(t)=\frac{1}{|\psi(P_j,x,y,t)|}$.

From the formula for the function $\psi$ follows that if the map
${\cal H}\circ \widetilde{r}$ is periodic for $t=0$, then it is
periodic with the same periods for any $t$ . The function
$\widetilde{r_j}$ satisfies the equations of Theorem 1 (see [1]).
Theorem 1 is proven.

We give an example of a Riemannian surface with the involutions $\sigma$ and $\tau$.
Let $\Gamma$ be a smooth supplement of the surface given in the $(\lambda,\mu)$-plane
by the equation
$$
\mu^3=\mu Q_1(\lambda)+Q_2(\lambda),
$$
where
$$
 Q_1(\lambda)=q_{-2k}\lambda^{-2k}+\dots+q_{2k}\lambda^{2k},\ \bar{q}_{-j}=q_j,
$$
$$
 Q_2(\lambda)=p_{-(2n+1)}\lambda^{-(2n+1)}+\dots+p_{2n+1}\lambda^{2n+1},\ \bar{p}_{-j}=-p_j.
$$
The surface $\Gamma$ has the holomorphic involution
$$
\sigma=(\lambda,\mu)=(-\lambda,-\mu)
$$
with two fixed points
$\infty_1=(0,\infty)$ and $\infty_2=(\infty,\infty)$
and the anti-holomorphic involution
$$
 \tau(\lambda,\mu)=\left(-\frac{1}{\bar{\lambda}},-\bar{\mu}\right).
$$

\newpage

{\bf Bibliography }

\vskip3mm

[1] Veselov A.P., Novikov, S.P. Finite-zone two-dimensional
potential Schr\"o\-dinger operators. Explicit formulas and
evolution equations // Sov. Math. Dokl. 1984, V. 30, P. 588--591.

[2] Taimanov I.A. The Weierstrass representation of closed
surfaces in ${\mathbb R}^3$ // Functional Anal. Appl. 1998, V. 32,
N. 4. P. 49--62.

[3] Dubrovin B.A., Krichever I.M., Novikov S.P.
The Schr\"odinger equation in a periodic field and Riemann surfaces //
Sov. Math. Dokl. 1976. V. 17. P. 947--951.

[4] Mironov A.E. New examples of Hamiltonian--minimal and minimal Lagrangian submanifolds in
${\mathbb C}^n$ and ${\mathbb C}P^n$. Sbornik Math. 2004. V. 195. N. 1. P. 85--96.

[5] Castro I., and Urbano F. Examples of unstable
Hamiltonian-minimal Lagrangian tori in ${\mathbb C}^2$ //
Compositio Math. 1998. V. 111. P. 1--14.

[6] Joyce D. Special Lagrangian 3-folds and integrable Systems // Advanced studies in pure math.
 Math. Soc. Japan (to appear).

[7] Sharipov R.A. Minimal tori in five-dimensional sphere in $C^3$ //
Theor. and Math. Phys. 1991. V. 87. N. 1. P. 48--56.

[8] Ma H., Ma J. Totally Real Minimal Tori in $CP^2$ // arXive: math.DG /0106141.

[9] Carberry E., McIntosh I. Minimal Lagrangian 2-tori in ${\mathbb C}P^2$ come
in real families of every dimension // London J. of Math. 2004. V. 69. N.2. P.531--544.

[10] Hitchin N. Harmonic maps from a 2-torus to the 3-spheres // J. diff. geom. 1990.
V. 31. P. 627--710.

[11] Mikhailov A.V. The reduction problem and the scattering method //Physica 3D. 1981. N. 1.
P. 73--117.

[12] Dubrovin B.A., Matveev V.B., Novikov, S.P. Non-linear
equations of Korteweg-de Vries type, finite-zone linear operators,
and Abelian varieties
// Russ. Math. Surv. 1976. V. 31, N.1, P. 59--146.

[13] Konopelchenko B.G. Induced surfaces and their integrable dynamics // Stud. Appl. Math. 1996. V. 96. N. 1. P. 9--51.

[14] Taimanov I.A. Modified Novikov-Veselov equation and differential geometry of surfaces //
Amer. Math. Soc. Transl. 1997. V. 179. Ser. 2. P. 131--151.

[15] Mironov A.E. On Hamiltonian--minimal Lagrangian tori in ${\mathbb C}P^2$ // Sib. Math. J.
2003. V. 44, N.6. P. 1039--1324.

[16] Krichever I.M. Nonlinear equations and elliptic curves //  J. Sov. Math. 1985. 28.
P. 51--90.
$$
$$
$$
$$

Andrey Mironov

Sobolev Institute of Mathematics SB RAS

Pr. acad. Koptyuga 4, 630090, Novosibirsk, Russia.

E-mail: mironov@math.nsc.ru

\end{document}